\documentclass{elsart}
\usepackage{amssymb,amsfonts}
\usepackage{graphicx}

\newcommand{\twomat}[4]{\left(\begin{array}{cc}#1&#2\\#3&#4\end{array}\right)}

\newcommand{\C}{{\mathbb{C}}}

\newcommand{\be}{\begin{equation}}
\newcommand{\ee}{\end{equation}}
\newcommand{\bea}{\begin{eqnarray}}
\newcommand{\eea}{\end{eqnarray}}
\newcommand{\beas}{\begin{eqnarray*}}
\newcommand{\eeas}{\end{eqnarray*}}

\newtheorem{theorem}{Theorem}

\newcount\minute
\newcount\hour
\def\currenttime{%
    \minute\time
    \hour\minute
    \divide\hour60
    \the\hour:\multiply\hour60\advance\minute-\hour\the\minute}
\begin{document}
\begin{frontmatter}
\title{Interpolating between the Arithmetic-Geometric Mean and Cauchy-Schwarz matrix norm inequalities}
\author{Koenraad M.R.\ Audenaert}
\address{
Department of Mathematics,
Royal Holloway University of London, \\
Egham TW20 0EX, United Kingdom \\[1mm]
Department of Physics and Astronomy, Ghent University, \\
S9, Krijgslaan 281, B-9000 Ghent, Belgium}
\ead{koenraad.audenaert@rhul.ac.uk}
\date{\today, \currenttime}
\begin{abstract}
We prove an inequality for unitarily invariant norms that interpolates between the Arithmetic-Geometric Mean
inequality and the Cauchy-Schwarz inequality.
\end{abstract}
\end{frontmatter}
\section{Introduction}
In this paper we prove the following inequality for unitarily invariant matrix norms:
\begin{theorem}\label{th:1}
Let $|||\cdot|||$ be any unitarily invariant norm.
For all $n\times n$ matrices $X$ and $Y$, and all $q\in[0,1]$,
\be
|||XY^*|||^2 \le |||qX^*X+(1-q)Y^*Y|||\;\;\;|||(1-q)X^*X+qY^*Y|||.\label{eq:1}
\ee
\end{theorem}
For $q=0$ or $q=1$, this reduces to the known Cauchy-Schwarz (CS) inequality for unitarily invariant norms
(\cite{bhatia}, inequality (IX.32))
\[
|||XY^*|||^2 \le |||X^*X|||\;\;\;|||Y^*Y|||.
\]
For $q=1/2$ on the other hand, this yields the arithmetic-geometric mean (AGM) inequality
(\cite{bhatia}, inequality (IX.22))
\[
|||XY^*||| \le \frac{1}{2}|||X^*X+Y^*Y|||.
\]
Thus, inequality (\ref{eq:1}) interpolates between the AGM and CS inequalities for unitarily invariant norms.

In Section 2 we prove an eigenvalue inequality that may be of independent interest. The proof of Theorem 1 follows easily from
this inequality, in combination with standard majorisation techniques; this proof is given in Section 3.
\section{Main technical result}
For any $n\times n$ matrix $A$ with real eigenvalues, we will denote these eigenvalues sorted in non-ascending order
by $\lambda_k(A)$. Thus $\lambda_1(A)\ge\cdots\ge\lambda_n(A)$.
Singular values will be denoted as $\sigma_k(A)$, again sorted in non-ascending order.

Our main technical tool in proving Theorem \ref{th:1} is the following eigenvalue inequality, 
which may be of independent interest:
\begin{theorem}\label{th:tool}
Let $A$ and $B$ be $n\times n$ positive semidefinite matrices.
Let $q$ be a number between 0 and 1, and let $C(q):=qA+(1-q)B$.
Then, for all $k=1,\ldots,n$,
\be
\lambda_k(AB) \le \lambda_k(C(q)C(1-q)).\label{eq:tool}
\ee
\end{theorem}
Putting $A=X^*X$ and $B=Y^*Y$, for $n\times n$ matrices $X$ and $Y$, and noting that
\[
\lambda_k^{1/2}(AB) = \lambda_k^{1/2}(YX^*XY^*)=\sigma_k(XY^*),
\]
we can write (\ref{eq:tool}) as a singular value inequality:
\be
\sigma_k^2(XY^*) \le \lambda_k((qX^*X +(1-q)Y^*Y)((1-q)X^*X+qY^*Y)).
\label{eq:tool2}
\ee

For $p=1/2$, Theorem \ref{th:tool} gives
\be
\lambda_k^{1/2}(AB) \le \frac{1}{2}\;\lambda_k(A+B)
\ee
and (\ref{eq:tool2}) becomes the well-known AGM inequality for singular values \cite[inequality (IX.20)]{bhatia}
\[
\sigma_k(XY^*) \le \frac{1}{2}\;\sigma_k(X^*X+Y^*Y).
\]

The following modification of inequality (\ref{eq:tool}) does not hold:
\[
\sigma_k(AB) \le \sigma_k(C(q)C(1-q)).
\]
We are grateful to Swapan Rana for finding counterexamples.

\bigskip

\noindent\textit{Proof of Theorem \ref{th:tool}.}
We first reduce the statement of the theorem to a special case
using a technique that is due to Ando \cite{ando95} and that was also used in \cite[Section 4]{drury12}.

Throughout the proof, we will keep $k$ fixed.
If either $A$ or $B$ has rank less than $k$, then $\lambda_k(AB)=0$ and (\ref{eq:tool}) holds trivially.
We will therefore assume that $A$ and $B$ have rank at least $k$.
By scaling $A$ and $B$ we can ensure that $\lambda_k(AB)=1$.

We will now try and find a positive semidefinite matrix $B'$ of rank exactly $k$ with $B'\le B$
and such that $AB'$ has $k$ eigenvalues equal to 1 and all others equal to 0.
By hypothesis, $AB$ and hence $A^{1/2} B A^{1/2}$ have at least $k$ eigenvalues larger than or equal to 1. 
Therefore, there exists a rank-$k$ projector $P$ satisfying $P\le A^{1/2} B A^{1/2}$.
Let $B'$ be a rank-$k$ matrix such that $A^{1/2} B' A^{1/2}=P$.
If $A$ is invertible, we simply have $B'=A^{-1/2} P A^{-1/2}$; otherwise the generalised inverse of $A$ is required.
Then $B'\le B$ and $AB'$ has the requested spectrum.

Passing to an eigenbasis of $B'$, we can decompose $B'$ as the direct sum
$B'=B_{11}\oplus[0]_{n-k}$, where $B_{11}$ is a $k\times k$ positive definite block.
In that same basis, we partition $A$ conformally with $B'$ as
$A=\twomat{A_{11}}{A_{12}}{A_{12}^*}{A_{22}}$.
Since $A^{1/2} B' A^{1/2}=P$ is a rank $k$ projector, so is 
\[
R:=(B')^{1/2} A (B')^{1/2}=(B_{11})^{1/2} A_{11} (B_{11})^{1/2}\oplus [0]_{n-k}.
\]
The top-left block of $R$ is a $k\times k$ matrix, and $R$
is a rank-$k$ projector. Therefore, that block must be identical to the $k\times k$ identity matrix:
$(B_{11})^{1/2} A_{11} (B_{11})^{1/2}=I$. This implies that $A_{11}$ is invertible and $B_{11}=(A_{11})^{-1}$.
We therefore have, in an eigenbasis of $B'$,
\[
A=\twomat{A_{11}}{A_{12}}{A_{12}^*}{A_{22}},\quad
B'=\twomat{(A_{11})^{-1}}{}{}{0}\le B.
\]

Clearly, $C'(q):= qA+(1-q)B'$ satisfies $C'(q)\le C(q)$, so that
\[
\lambda_k(C'(q)C'(1-q))\le \lambda_k(C(q)C(1-q)),
\] 
while still $\lambda_k(AB')=\lambda_k(AB)=1$.
It is now left to show that $\lambda_k(C'(q)C'(1-q))\ge 1$.

A further reduction is possible.
Let
\[
A' = \twomat{A_{11}}{A_{12}}{A_{12}^*}{\quad A_{12}^*(A_{11})^{-1}A_{12}},
\]
which has rank $k$ and satisfies $0\le A'\le A$.
Let also $C''(q):=qA'+(1-q)B'$, for which $0\le C''(q)\le C'(q)$.
Then $\lambda_k(C''(q)C''(1-q))\le \lambda_k(C'(q)C'(1-q))$.

Introducing $F:=A_{11}>0$, $G:=A_{12}A_{12}^*\ge0$ and $s:=(1-q)/q>0$,
we have
\beas
C''(q) &=& q\twomat{F}{A_{12}}{A_{12}^*}{A_{12}^* F^{-1}A_{12}} + (1-q) \twomat{F^{-1}}{}{}{0} \\
&=& q\twomat{I}{}{}{A_{12}^*}\;\twomat{F+sF^{-1}}{I}{I}{F^{-1}} \;\twomat{I}{}{}{A_{12}}
\eeas
so that
\beas
\lefteqn{\lambda_k(C''(q)C''(1-q))} \\ 
&=& q(1-q)\lambda_k\left(
\twomat{I}{}{}{G}
\twomat{F+sF^{-1}}{I}{I}{F^{-1}}
\twomat{I}{}{}{G}
\twomat{F+s^{-1}F^{-1}}{I}{I}{F^{-1}}
\right),
\eeas
where each factor is a $2k\times 2k$ matrix.
Noting that
\[
\twomat{F+sF^{-1}}{I}{I}{F^{-1}}
= \twomat{s^{1/2} F^{-1/2}}{F^{1/2}}{0}{F^{-1/2}} \;
\twomat{s^{1/2} F^{-1/2}}{0}{F^{1/2}}{F^{-1/2}},
\]
we then have
$\lambda_k(C''(q)C''(1-q)) = q(1-q)\lambda_k(Z^*Z) = q(1-q)\sigma_k^2(Z)$, where
\[
Z= \twomat{s^{1/2} F^{-1/2}}{0}{F^{1/2}}{F^{-1/2}}\;
\twomat{I}{}{}{G}\;
\twomat{s^{-1/2}F^{-1/2}}{F^{1/2}}{0}{F^{-1/2}}
= \twomat{F^{-1}}{s^{1/2}}{s^{-1/2}}{F+H},
\]
and $H:=F^{-1/2}GF^{-1/2}\ge0$.
The singular values of $Z$ are the same as those of 
\[
X:= \twomat{s^{1/2}}{F^{-1}}{F+H}{s^{-1/2}}.
\]
We now use the fact that the singular values of $X$ are bounded below by the ordered eigenvalues
of the Hermitian part of $X$: $\sigma_j(X)\ge\lambda_j((X+X^*)/2)$ for $j=1,\ldots,2k$ \cite[Corollary 3.1.5]{HJII}.
Thus,
\[
\lambda_k(C''(q)C''(1-q)) \ge q(1-q)\lambda_k^2(Y),
\]
\[
\mbox{ with } Y:=\twomat{s^{1/2}}{K}{K}{s^{-1/2}}
\mbox{ and } K:=(F+H+F^{-1})/2.
\]
Clearly, $K\ge (F+F^{-1})/2 \ge I$.
It is easily checked that the $k$ largest eigenvalues of $Y$
are given by
\[
\lambda_j(Y)=\frac{1}{2}\left(s^{1/2}+s^{-1/2}+\sqrt{(s^{1/2}+s^{-1/2})^2-4+4\lambda_j^2(K)}\right),\quad j=1,\ldots,k.
\]
As this expression is a monotonously increasing function of $\lambda_j(K)$, and $\lambda_j(K)\ge1$, we obtain the lower bound
$\lambda_k(Y)\ge s^{1/2}+s^{-1/2}$.
Then, finally,
\beas
\lambda_k(C''(q)C''(1-q))
&\ge& q(1-q)\;(s^{1/2}+s^{-1/2})^2 \\
&=& q(1-q)\left(\left(\frac{1-q}{q}\right)^{1/2}+\left(\frac{q}{1-q}\right)^{1/2}\right)^2 \\
&=& (1-q+q)^2 = 1,
\eeas
from which it follows that $\lambda_k(C'(q)C'(1-q))\ge 1$.
\qed
\section{Proof of Theorem \ref{th:1}}
Using Theorem \ref{th:tool} and some standard arguments, the promised norm inequality is easily proven.

For all positive semidefinite matrices $A$ and $B$, and any $r>0$, we have the weak majorisation relation
\[
\lambda^r(AB) \prec_w \lambda^r(A)\cdot\lambda^r(B),
\]
where `$\cdot$' denotes the elementwise product for vectors.
This relation follows from combining the fact that $AB$ has non-negative eigenvalues with Weyl's majorant inequality
(\cite{bhatia}, inequality (II.23))
\[
|\lambda(AB)|^r \prec_w \sigma^r(AB)
\]
and with the singular value majorisation relation
(\cite{bhatia}, inequality (IV.41))
\[
\sigma^r(AB) \prec_w \sigma^r(A)\cdot\sigma^r(B).
\]

From (\ref{eq:tool2}) we immediately get, for any $r>0$,
\[
\sigma^{2r}(XY^*) \prec_w \lambda^r\left((qX^*X +(1-q)Y^*Y)\;((1-q)X^*X+qY^*Y)\right).
\]
Hence,
\[
\sigma^{2r}(XY^*) \prec_w \lambda^r(qX^*X +(1-q)Y^*Y)\cdot\lambda^r((1-q)X^*X+qY^*Y)).
\]

If we now apply H\"older's inequality for symmetric gauge functions $\Phi$,
\[
\Phi(|x\cdot y|) \le \Phi(|x|^p)^{1/p}\;\Phi(|y|^{p'})^{1/p'},
\]
where $x,y\in\C^n$ and $1/p+1/p'=1$,
we obtain
\beas
\Phi( \sigma^{2r}(XY^*) ) &\le& \Phi( \lambda^r(qX^*X +(1-q)Y^*Y)\cdot \lambda^r((1-q)X^*X+qY^*Y)) ) \\
&\le& \Phi( \lambda^{rp}(qX^*X +(1-q)Y^*Y))^{1/p}\;\Phi( \lambda^{rp'}((1-q)X^*X+qY^*Y)) )^{1/p'}.
\eeas
Hence, for any unitarily invariant norm,
\[
||| \;|XY^*|^{2r} \;|||
\le |||(qX^*X +(1-q)Y^*Y)^{rp}|||^{1/p}\;|||((1-q)X^*X+qY^*Y)^{rp'}|||^{1/p'}.
\]
Theorem \ref{th:1} now follows by setting $r=1/2$ and $p=p'=2$. \qed
\section*{Acknowledgments}
We acknowledge support by an Odysseus grant from the Flemish FWO.
We are grateful to Professor Bhatia for pointing out a serious mistake in an earlier circulated version of this paper.

\end{document}